\documentclass[a4paper,notitlepage,11pt]{article}%
\usepackage{mitpress}
\usepackage[centertags]{amsmath}
\usepackage{graphicx}
\usepackage{amsfonts}
\usepackage{amssymb}%
\providecommand{\U}[1]{\protect\rule{.1in}{.1in}}
\providecommand{\U}[1]{\protect\rule{.1in}{.1in}}
\providecommand{\U}[1]{\protect\rule{.1in}{.1in}}
\providecommand{\U}[1]{\protect\rule{.1in}{.1in}}
\newtheorem{theorem}{Theorem}
\newtheorem{acknowledgement}[theorem]{Acknowledgement}

\newtheorem{corollary}[theorem]{Corollary}

\newtheorem{proposition}[theorem]{Proposition}
\newtheorem{remark}[theorem]{Remark}

\newenvironment{proof}[1][Proof]{\textbf{#1.} }{\ \rule{0.5em}{0.5em}}
\newdimen\dummy
\dummy=\oddsidemargin
\begin{document}

\title{Discrete Hybrid Polyharmonic Cubature Formulas with weight on the disc. with
error bounds}
\author{O. Kounchev and H. Render}
\maketitle

\begin{abstract}
This paper is a second part of our study of the Discrete Polyharmonic Cubature
Formulas on the disc. It completes our study and provides a satisfactory
cubature formula in terms of precision and number of evaluation points
(coefficient of efficiency of the formula). 

\end{abstract}

\section{ Introduction}

The present paper is a continuation of paper
\cite{kounchevRender2015DiscretePolyharmonicCubature}, where we have provided
a Discrete Polyharmonic Cubature Formula for integration of functions
multiplied by a weight on the disc. The main merit of this formula is that
unlike majority of research in the area, we provide constructive error bounds.
However an issue arises with the number of evaluation points. In the present
paper we show how to solve this issue by using spline approximation methods.
This approximation increases the error of the evaluation of the integral but
it is again under control and we show that it is a reasonable compromise which
gives much better result compared to the alternative methods.

We decided to provide a separate presentation of the spline version of the
Discrete Polyharmonic Cubature Formula, which we call here \textbf{Hybrid
Discrete Polyharmonic Cubature Formula}, in order not to overburden the
attention of the reader of the original paper
\cite{kounchevRender2015DiscretePolyharmonicCubature}. We will rely upon the
notations and definitions provided in
\cite{kounchevRender2015DiscretePolyharmonicCubature}.

As in \cite{kounchevRender2015DiscretePolyharmonicCubature}, we consider
numerical evaluation of integrals of the type%
\begin{equation}
I_{w}\left(  f\right)  =\int_{D_{R}}f\left(  x\right)  w\left(  x\right)
dx\label{eqint}%
\end{equation}
where $w\left(  x\right)  $ is a (not necessarily non-negative) weight
function on the disc $D_{R}$ in the plane $\mathbb{R}^{2}.$ As in
\cite{kounchevRender2015DiscretePolyharmonicCubature} we provide experiments
with the weight functions
\begin{equation}
w^{\left(  1\right)  }\left(  x,y\right)  =\frac{1+x}{\sqrt{x^{2}+y^{2}}}%
\quad\text{ and }\quad w^{\left(  2\right)  }\left(  x,y\right)  =\left\vert
y\right\vert \ .\label{eqweight}%
\end{equation}
We provide the experiments in a form facilitating comparison with the results
in \cite{kounchevRender2015DiscretePolyharmonicCubature}.

For completeness of notations in the present paper, we repeat the main
formulas used in the previous paper
\cite{kounchevRender2015DiscretePolyharmonicCubature}:

\begin{enumerate}
\item The Fourier series of $w$ is given by
\begin{equation}
w\left(  re^{i\varphi}\right)  =w\left(  r\cos\varphi,r\sin\varphi\right)  :=%
{\displaystyle\sum_{k=0}^{\infty}}
{\displaystyle\sum_{\ell=1}^{a_{k}}}
w_{\left(  k,\ell\right)  }\left(  r\right)  Y_{\left(  k,\ell\right)
}\left(  \varphi\right)  . \label{FourierSeries}%
\end{equation}

\item Orthonormalized spherical harmonics are defined by
\begin{align}
Y_{\left(  0,1\right)  }\left(  \varphi\right)   &  =1/\sqrt{2\pi
}\label{sphericalHarmonics1}\\
Y_{\left(  k,1\right)  }\left(  \varphi\right)   &  =\frac{1}{\sqrt{\pi}}\cos
k\varphi\quad\text{ and }\quad Y_{\left(  k,2\right)  }\left(  \varphi\right)
=\frac{1}{\sqrt{\pi}}\sin k\varphi. \label{sphericalHarmonics2}%
\end{align}
for integers $k\geq1.$ Then $Y_{\left(  k,\ell\right)  }$ is an
\emph{orthonormal} system for $k\geq0,\ \ell=1,..,a_{k}$, where $a_{k}=2$ for
$k\geq1$, and $a_{0}=1.$

\item The $\left(  k,\ell\right)  $-th Fourier coefficient of a complex-valued
continuous function $f\left(  re^{i\varphi}\right)  $ is given by
\begin{equation}
f_{\left(  k,l\right)  }\left(  r\right)  :=\int_{0}^{2\pi}f\left(
re^{i\varphi}\right)  Y_{\left(  k,\ell\right)  }\left(  \varphi\right)
d\varphi\qquad\text{ for }k\geq0,\ell=1,..,a_{k} \label{fkl}%
\end{equation}

\item The \emph{corresponding Fourier series} of $f$ is
\begin{equation}
f\left(  re^{i\varphi}\right)  =f\left(  r\cos\varphi,r\sin\varphi\right)  :=%
{\displaystyle\sum_{k=0}^{\infty}}
{\displaystyle\sum_{\ell=1}^{a_{k}}}
f_{\left(  k,\ell\right)  }\left(  r\right)  Y_{\left(  k,\ell\right)
}\left(  \varphi\right)  . \label{FourierSeriesGeneral}%
\end{equation}

\item The Fourier series of a polynomial $P\left(  x\right)  $ is of a very
special form: there exist polynomials $\widetilde{p}_{\left(  k,\ell\right)
}$ and a number $N\leq\deg P\left(  x\right)  $ such that
\begin{equation}
P\left(  x\right)  =P\left(  r\cos\varphi,r\sin\varphi\right)  =%
{\displaystyle\sum_{k=0}^{N}}
{\displaystyle\sum_{\ell=1}^{a_{k}}}
\widetilde{p}_{\left(  k,\ell\right)  }\left(  r^{2}\right)  r^{k}Y_{\left(
k,\ell\right)  }\left(  \varphi\right)  ;\label{Almansi}%
\end{equation}
the representation (\ref{Almansi}) is called Gauss decomposition or
\textbf{Almansi expansion} of a polynomial $p$.

\item Hence, the Fourier coefficient $p_{\left(  k,\ell\right)  }\left(
r\right)  $ of a polynomial $P\left(  x\right)  $ is of the form
\begin{equation}
p_{\left(  k,\ell\right)  }\left(  r\right)  =\widetilde{p}_{\left(
k,\ell\right)  }\left(  r^{2}\right)  r^{k}.\label{Almansi2}%
\end{equation}

\item The integral (\ref{eqint}), in polar coordinates, becomes%
\[
I_{w}\left(  f\right)  =\int_{0}^{2\pi}\int_{0}^{R}f\left(  r\cos\varphi
,r\sin\varphi\right)  \cdot w\left(  r\cos\varphi,r\sin\varphi\right)  \cdot
rdrd\varphi.
\]

\item We obtain
\begin{equation}
I_{w}\left(  f\right)  =\int_{0}^{2\pi}\int_{0}^{R}f\left(  re^{i\varphi
}\right)  w\left(  re^{i\varphi}\right)  rdrd\varphi=%
{\displaystyle\sum_{k=0}^{\infty}}
{\displaystyle\sum_{\ell=1}^{a_{k}}}
\int_{0}^{R}f_{k,\ell}\left(  r\right)  w_{\left(  k,\ell\right)  }\left(
r\right)  rdr.\label{eqcentral2}%
\end{equation}

\item For the constant weight function $w\left(  x,y\right)  =1=\sqrt{2\pi
}Y_{\left(  0,1\right)  }$ one obtains simply
\begin{equation}
I_{1}\left(  f\right)  =\int_{0}^{2\pi}\int_{0}^{R}f\left(  r\cos\varphi
,r\sin\varphi\right)  rdrd\varphi=\sqrt{2\pi}\int_{0}^{R}f_{\left(
0,1\right)  }\left(  r\right)  rdr.\label{eqintegralpolar}%
\end{equation}

\item We shall \textbf{assume} that each function $w_{\left(  k,\ell\right)
}\left(  r\right)  $\ does not change the sign over the interval $\left[
0,R\right]  ,$\ a property which we call \textbf{pseudo-definiteness}%
\begin{align}
w_{\left(  k,\ell\right)  }\left(  r\right)    & \geq0\qquad\text{on }\left[
0,R\right]  \text{ or }\label{ConditionPseudoDef}\\
w_{\left(  k,\ell\right)  }\left(  r\right)    & \leq0\qquad\text{on }\left[
0,R\right]  .\nonumber
\end{align}

\item By the Almansi formula (\ref{Almansi}), (\ref{Almansi2}), and a change
of the variable $\rho=r^{2},$  the one-dimensional integrals in
(\ref{eqcentral2}) used for computing $I_{w}\left(  p\right)  $ are equal to
\begin{equation}
\int_{0}^{R}p_{k,\ell}\left(  r\right)  w_{\left(  k,\ell\right)  }\left(
r\right)  rdr=\int_{0}^{R}\widetilde{p}_{k,\ell}\left(  r^{2}\right)
r^{k}w_{\left(  k,\ell\right)  }\left(  r\right)  rdr=\frac{1}{2}\int
_{0}^{R^{2}}P\left(  \rho\right)  \cdot\rho^{k/2}w_{\left(  k,\ell\right)
}\left(  \sqrt{\rho}\right)  d\rho.\label{defTkl}%
\end{equation}

\item Assuming \ref{ConditionPseudoDef} to the  integral (\ref{defTkl}) with
measure $\rho^{k/2}w_{\left(  k,\ell\right)  }\left(  \sqrt{\rho}\right)  ,$
we apply the $N$-point Gauss-Jacobi quadrature: we obtain the nodes
$t_{1,\left(  k,\ell\right)  }<...<t_{N,\left(  k,\ell\right)  }$ and the
weights $\lambda_{1,\left(  k,\ell\right)  },...,\lambda_{N,\left(
k,\ell\right)  }$ which are either all positive or all negative. Due to the
exactness of the Gauss-Jacobi quadrature for any integer $0\leq s\leq2N-1$ we
obtain the equalities
\begin{equation}%
{\displaystyle\sum_{j=1}^{N}}
\lambda_{j,\left(  k,\ell\right)  }\cdot t_{j,\left(  k,\ell\right)  }%
^{s}=\frac{1}{2}\int_{0}^{R^{2}}\rho^{s}\rho^{k/2}w_{\left(  k,\ell\right)
}\left(  \sqrt{\rho}\right)  d\rho\qquad\left(  =\int_{0}^{R}r^{2s}%
r^{k}w_{\left(  k,\ell\right)  }\left(  r\right)  rdr\right)
\label{eqexactkl}%
\end{equation}

\item Hence, for a polynomial $f,$ for which $\deg f_{\left(  k,\ell\right)
}\leq2N-1$ for all $\left(  k,\ell\right)  $, by using the Gauss-Jacobi
quadrature (\ref{eqexactkl}), the integral (\ref{eqcentral2}) becomes
\begin{align}
I_{w}\left(  f\right)   &  =%
{\displaystyle\sum_{k=0}^{\infty}}
{\displaystyle\sum_{\ell=1}^{a_{k}}}
\int_{0}^{R}f_{k,\ell}\left(  r\right)  w_{\left(  k,\ell\right)  }\left(
r\right)  rdr\label{Iwf2}\\
&  =\frac{1}{2}%
{\displaystyle\sum_{k=0}^{\infty}}
{\displaystyle\sum_{\ell=1}^{a_{k}}}
\int_{0}^{R^{2}}f_{k,\ell}\left(  \sqrt{\rho}\right)  \rho^{-k/2}\left\{
\rho^{k/2}w_{\left(  k,\ell\right)  }\left(  \sqrt{\rho}\right)  \right\}
d\rho\nonumber\\
&  =I_{N}^{\text{poly}}\left(  f\right)  .\nonumber
\end{align}
where we have put
\begin{equation}
I_{N}^{\text{poly}}\left(  f\right)  :=\frac{1}{2}%
{\displaystyle\sum_{k=0}^{\infty}}
{\displaystyle\sum_{\ell=1}^{a_{k}}}
{\displaystyle\sum_{j=1}^{N}}
\lambda_{j,\left(  k,\ell\right)  }\cdot t_{j,\left(  k,\ell\right)  }%
^{-\frac{k}{2}}\cdot f_{\left(  k,\ell\right)  }\left(  \sqrt{t_{j,\left(
k,\ell\right)  }}\right)  ;\label{eqIpoly}%
\end{equation}

\item As the values $f_{k,\ell}\left(  \sqrt{t_{j,\left(  k,\ell\right)  }%
}\right)  $ are Fourier coefficients, they can be approximated by means of
Discrete Fourier transform of the function $f.$ Thus for every two integers
$K,M\geq0$ we introduced the \textbf{Discrete Polyharmonic Cubature with
parameters} $\left(  N,M,K\right)  $ in the paper
\cite{kounchevRender2015DiscretePolyharmonicCubature}, by putting
\begin{align}
I_{\left(  N,M,K\right)  }^{\text{poly}}\left(  f\right)   &  :=\frac{1}{2}%
{\displaystyle\sum_{k=0}^{K}}
{\displaystyle\sum_{\ell=1}^{a_{k}}}
{\displaystyle\sum_{j=1}^{N}}
\lambda_{j,\left(  k,\ell\right)  }\cdot t_{j,\left(  k,\ell\right)  }%
^{-\frac{k}{2}}\cdot f_{\left(  k,\ell\right)  }^{\left(  M\right)  }\left(
\sqrt{t_{j,\left(  k,\ell\right)  }}\right)  \label{DiscrPolyhCubature}\\
&  =\frac{\pi}{M}%
{\displaystyle\sum_{k=0}^{K}}
{\displaystyle\sum_{\ell=1}^{a_{k}}}
{\displaystyle\sum_{j=1}^{N}}
{\displaystyle\sum_{s=1}^{M}}
\lambda_{j,\left(  k,\ell\right)  }\cdot t_{j,\left(  k,\ell\right)  }%
^{-\frac{k}{2}}\cdot Y_{\left(  k,\ell\right)  }\left(  \frac{2\pi s}%
{M}\right)  \cdot f\left(  \sqrt{t_{j,\left(  k,\ell\right)  }}e^{i\frac{2\pi
s}{M}}\right)  \label{DiscrPolyhCubature+}%
\end{align}
where the Discrete Fourier transform is given by
\begin{equation}
f_{\left(  k,\ell\right)  }^{\left(  M\right)  }\left(  r\right)  :=\frac
{2\pi}{M}%
{\displaystyle\sum_{s=1}^{M}}
f\left(  re^{i\frac{2\pi s}{M}}\right)  Y_{\left(  k,\ell\right)  }\left(
\frac{2\pi s}{M}\right)  \label{DFT}%
\end{equation}

\item Its coefficients (weights) $\left\{  \lambda_{j,\left(  k,\ell\right)
}\cdot t_{j,\left(  k,\ell\right)  }^{-\frac{k}{2}}\cdot Y_{\left(
k,\ell\right)  }\left(  \frac{2\pi s}{M}\right)  \right\}  $ have varying
signs but they satisfy the following \textbf{remarkable} inequality
\begin{equation}
\frac{\pi}{M}%
{\displaystyle\sum_{k=0}^{K}}
{\displaystyle\sum_{\ell=1}^{a_{k}}}
{\displaystyle\sum_{j=1}^{N}}
{\displaystyle\sum_{s=1}^{M}}
\left\vert \lambda_{j,\left(  k,\ell\right)  }\cdot t_{j,\left(
k,\ell\right)  }^{-\frac{k}{2}}\cdot Y_{\left(  k,\ell\right)  }\left(
\frac{2\pi s}{M}\right)  \right\vert \leq\sqrt{\pi}\left\Vert w\right\Vert
,\label{wRemarkable}%
\end{equation}
where we have assumed that the weight $w$ satisfies
\begin{equation}
\left\Vert w\right\Vert :=\sum_{k=0}^{\infty}\sum_{\ell=1}^{a_{k}}\int_{0}%
^{R}\left\vert w_{\left(  k,\ell\right)  }\left(  r\right)  \right\vert
rdr<\infty.\label{wNorm}%
\end{equation}
This inequality is proved in
\cite{kounchevRender2015DiscretePolyharmonicCubature}, by an application of
the famous Chebyshev extremal property for the Gauss-Jacobi quadrature:
\begin{equation}%
{\displaystyle\sum_{j=1}^{N}}
\left\vert \lambda_{j,\left(  k,\ell\right)  }\right\vert t_{j,\left(
k,\ell\right)  }^{-\frac{k}{2}}\leq\frac{1}{2}\int_{0}^{R^{2}}\left\vert
w_{\left(  k,\ell\right)  }\left(  \sqrt{\rho}\right)  \right\vert d\rho
=\int_{0}^{R}\left\vert w_{\left(  k,\ell\right)  }\left(  r\right)
\right\vert rdr.\label{GNkl}%
\end{equation}

\item Let us note that formula (\ref{DiscrPolyhCubature+}) needs $\left(
2K-1\right)  \cdot N\cdot M$ evaluation points. The coefficient of efficiency
of a cubature formula may be introduced in the following way (see e.g.
\cite{xiaoGimbutas}): If $n$ is the number of nodes of the quadrature, $m$ is
the number of linearly independent functions which are integrated exactly, and
$d$ is the dimension of the Euclidean space, then the coefficient $E$ is
defined as%
\[
E=\frac{m}{n\left(  d+1\right)  }.
\]
For Gaussian quadratures in one dimension, $E$ is obviously equal to $1$. For
the Discrete Polyharmonic Cubature formula (\ref{DiscrPolyhCubature}) we need
precisely $\left(  2K-1\right)  \cdot N\cdot M$ sampling points -- hence $2$
coordinates for every point and $1$ coefficient. On the other hand, as we have
seen in \cite{kounchevRender2015DiscretePolyharmonicCubature}, the subspace of
exactness of formula (\ref{DiscrPolyhCubature}) has dimension $2N\times\left(
2M-3-2K\right)  $, hence the efficiency coefficient of the Discrete
Polyharmonic Cubature formula is
\[
E=\frac{2N\times\left(  2M-3-2K\right)  }{3\left(  2K-1\right)  \times N\times
M}%
\]
which is much less than $1,$ hence is prettily bad. The main purpose of the
present paper is to show how by using spline approximation and sacrificing
some precision, it is possible to reduce the number of point evaluations which
will improve essentially the coefficient $E,$ by making it very close to $1.$
\end{enumerate}

Now we come to the main point of the present paper.

Since the point evaluations of a function $f$ can be very costly, it is of
advantage to reduce their number by sacrificing some precision. There are
simple ways to decrease essentially the number of evaluation points, e.g. by
applying spline approximation to the coefficients $f_{\left(  k,\ell\right)
}\left(  \sqrt{t_{j,\left(  k,\ell\right)  }}\right)  .$ For example, let us
assume that for some integer $N_{1}\geq1,$ we are given the radii $R_{j},$
$j=1,2,...,N_{1}$ with $0<R_{1}<R_{2}<\cdot\cdot\cdot<R_{N_{1}}\leq R,$ and
the $N_{1}\times M$ functional values be given as well:
\[
f\left(  R_{j}e^{i\varphi_{s}}\right)  \qquad\text{for }%
s=1,2,...,M,\ \text{and }j=1,2,...,N_{1}.
\]
Thus we will have only $N_{1}\times M$ sampling points of the function $f$
which is the usual framework for the cubature formulas. Then, by formula
(\ref{DFT}), we will find the Discrete Fourier transform approximation
$f_{k,\ell}^{\left(  M\right)  }\left(  R_{j}\right)  $ to $f_{k,\ell}\left(
R_{j}\right)  ,$ and we may use \textbf{interpolation splines} with data
$\left\{  f_{k,\ell}\left(  R_{j}\right)  \right\}  _{j=1}^{N_{1}}$ to
approximate all values $f_{k,\ell}^{\left(  M\right)  }\left(  \sqrt
{t_{j,\left(  k,\ell\right)  }}\right)  $ for $j=1,2,...,N.$ The resulting
formula introduced below is termed \textbf{(Spline) Hybrid  Discrete
Polyharmonic Cubature} formula (\ref{hybridCubature}).

As with the experiments in
\cite{kounchevRender2015DiscretePolyharmonicCubature}, the experiments with
formula (\ref{hybridCubature}), show excellent performance. This is due to the
fact that we are able to control the error bounds of formula
(\ref{hybridCubature}), where we now have the parameter $N_{1},$ for the knot
evaluations, see Theorem \ref{ThybridError} below.

The plan of the paper is the following: In Section \ref{Shybrid} we define the
spline based Hybrid Discrete Polyharmonic Cubature and provide error bounds
for it.

In Section \ref{Sexperiments} we provide experimental results for the discrete
polyharmonic cubature with respect to the first weight function $w^{\left(
1\right)  }$ in (\ref{eqweight}).

Section \ref{Sexperiments2} contains experimental results for the second
weight function $w^{\left(  2\right)  }.$

\section{The Hybrid Polyharmonic Cubature formula \label{Shybrid}}

Let let us assume that on the two-dimensional disc $B_{R}$ of radius $R,$  the
values of the function $f\left(  x\right)  ,$ $x\in B_{R}\subset\mathbb{R}%
^{2}$  are given on a regular grid: We assume that the integers $N_{1}$ and
$M$ are fixed in advance, where $M$ is odd. Then we consider the  concentric
circles with radii $R_{j},$ for $j=1,2,...,$ $N_{1}$, with $0=R_{0}%
<R_{1}<R_{2}<\cdot\cdot\cdot<R_{N_{1}}=R,$ and we assume that the values of
the function $f$ are given:\
\[
f\left(  R_{m}e^{i\frac{2\pi}{M}s}\right)  \qquad\text{for }m=1,2,...,N_{1}%
,\quad s=1,2,...,M.
\]
Formula (\ref{DFT}) provides an approximation $f_{k,\ell}^{\left(  M\right)
}\left(  R_{m}\right)  $ to the functions $f_{k,\ell}\left(  R_{m}\right)  .$

The value $f_{k,\ell}\left(  r\right)  $ at the knot $r=0$ is determined from
the following arguments: In the case of $k\geq1,$ by Proposition $1$ in
\cite{baouendiGoulaouicLipkin}, we know that if the function $f\in C^{\infty
}\left(  \overline{D_{R}}\right)  $ then $f_{\left(  k,\ell\right)  }\in
C^{\infty}\left[  0,R\right]  ,$ and
\[
f_{\left(  k,\ell\right)  }\left(  0\right)  =0.
\]
Hence, we will put also
\[
f_{\left(  k,\ell\right)  }^{\left(  M\right)  }\left(  0\right)
=0\qquad\text{for }k\geq1.
\]
On the other hand, we see by formulas (\ref{sphericalHarmonics1}), (\ref{fkl})
that the very first coefficient satisfies
\[
f_{\left(  0,1\right)  }\left(  0\right)  =\frac{1}{\sqrt{2\pi}}f\left(
0\right)  ,
\]
hence, we will have
\[
f_{\left(  0,1\right)  }^{\left(  M\right)  }\left(  0\right)  =\frac{1}%
{\sqrt{2\pi}}f\left(  0\right)  .
\]

In general, if $V\in\mathbb{C}^{N_{1}+1}$ is a vector, we will denote by
\begin{equation}
SPL\left[  V\right]  \left(  r\right)  \label{spline}%
\end{equation}
the value at $r,$ $0\leq r\leq1,$ of the interpolation cubic spline with knots
$\left\{  R_{j}\right\}  _{j=0}^{N_{1}}$ and interpolation data $\left\{
V_{j}\right\}  _{j=0}^{N_{1}}$, which satisfies the \textbf{not-a-knot}
boundary conditions, see Proposition \ref{PBeatsonDeBoor} below. If we denote
by $S_{m}\left(  r\right)  $ the\emph{ not-a-knot} interpolation spline which
satisfies
\[
S_{m}\left(  R_{s}\right)  =\delta_{m,s}%
\]
(where $\delta$ denotes the Kronecker symbol) then we obtain the
representation for the data $V,$ given by
\[
SPL\left[  V\right]  \left(  r\right)  =\sum_{m=0}^{N_{1}}S_{m}\left(
r\right)  V_{m}.
\]

The idea is now to approximate the values $f_{k,\ell}^{\left(  M\right)
}\left(  \sqrt{t_{j,\left(  k,\ell\right)  }}\right)  $ through the
\emph{not-a-knot} spline $SPL$ interpolating the values $\left\{  f_{\left(
k,\ell\right)  }^{\left(  M\right)  }\left(  R_{j}\right)  \right\}
_{j=0}^{N_{1}}$ at the nodes $0=R_{0}<R_{1}<R_{2}<\cdot\cdot\cdot<R_{N_{1}%
}=R.$ Thus, from the Discrete Polyharmonic Cubature with parameters $\left(
N,M,K\right)  $ defined in (\ref{DiscrPolyhCubature}), after substituting
$f_{\left(  k,\ell\right)  }^{\left(  M\right)  }\left(  \sqrt{t_{j,\left(
k,\ell\right)  }}\right)  $ with $SPL\left[  \left\{  f_{\left(
k,\ell\right)  }^{\left(  M\right)  }\left(  R_{m}\right)  \right\}
_{m=0}^{N_{1}}\right]  \left(  \sqrt{t_{j,\left(  k,\ell\right)  }}\right)  $
we obtain the \textbf{Hybrid Polyharmonic Cubature} formula:
\begin{equation}
I_{\left(  N,M,K,N_{1}\right)  }^{\text{spline}}\left(  f\right)  :=\frac
{1}{2}%
{\displaystyle\sum_{k=0}^{K}}
{\displaystyle\sum_{\ell=1}^{a_{k}}}
{\displaystyle\sum_{j=1}^{N}}
\lambda_{j,\left(  k,\ell\right)  }\times t_{j,\left(  k,\ell\right)
}^{-\frac{k}{2}}\times SPL\left[  \left\{  f_{\left(  k,\ell\right)
}^{\left(  M\right)  }\left(  R_{m}\right)  \right\}  _{m=0}^{N_{1}}\right]
\left(  \sqrt{t_{j,\left(  k,\ell\right)  }}\right)  \label{hybridCubature}%
\end{equation}
Obviously, we have the equalities:\
\begin{align}
I_{\left(  N,M,K,N_{1}\right)  }^{\text{spline}}\left(  f\right)   &
=\frac{1}{2}\sum_{k=0}^{K}\sum_{\ell=1}^{a_{k}}\sum_{j=1}^{N}\sum_{m=0}%
^{N_{1}}t_{j,\left(  k,\ell\right)  }^{-\frac{k}{2}}\lambda_{j,\left(
k,\ell\right)  }\times S_{m}\left(  \sqrt{t_{j,\left(  k,\ell\right)  }%
}\right)  f_{\left(  k,\ell\right)  }^{\left(  M\right)  }\left(  R_{m}\right)
\label{hybridCubature+}\\
&  =\sum_{m=0}^{N_{1}}\sum_{s=1}^{M}f\left(  R_{m}e^{i\frac{2\pi}{M}s}\right)
w_{m,s}\nonumber
\end{align}
where we have put
\[
w_{m,s}:=\frac{1}{2}\sum_{k=0}^{K}\sum_{\ell=1}^{a_{k}}\sum_{j=1}^{N}%
S_{m}\left(  \sqrt{t_{j,\left(  k,\ell\right)  }}\right)  \frac{2\pi}%
{M}Y_{\left(  k,\ell\right)  }\left(  \frac{2\pi}{M}s\right)  \times
t_{j,\left(  k,\ell\right)  }^{-\frac{k}{2}}\lambda_{j,\left(  k,\ell\right)
}.
\]
The last formula (\ref{hybridCubature+}) shows that $I_{\left(  N,M,K,N_{1}%
\right)  }^{\text{spline}}\left(  f\right)  $ is a cubature formula in the
usual sense, having $N_{1}\cdot M$ nodes.

Let us note that formula (\ref{hybridCubature+}) has the same knots as the
Midpoint cubature rule, see
\cite{kounchevRender2015DiscretePolyharmonicCubature}.

Let us put
\begin{equation}
h:=\max_{i}\left(  R_{i+1}-R_{i}\right)  . \label{hR}%
\end{equation}

Then the error bound for not-a-knot spline is given in the following
Proposition. 

\begin{proposition}
\label{PBeatsonDeBoor}Let $g\in C^{4}\left(  \left[  a,b\right]  \right)  $.
Assume that the knots of the interpolation are $a=R_{0}<R_{1}<\cdot\cdot
\cdot<R_{N_{1}}=b.$ Let $N_{1}\geq5.$ Let $s\left(  t\right)  $ be the
not-a-knot intepolation cubic spline, i.e. the $C^{2}$ piece-wise cubic
polynomial function which satisfies the interpolation conditions
\[
s\left(  R_{j}\right)  =g\left(  R_{j}\right)  \qquad\text{for }%
j=0,1,...,N_{1}%
\]
and the \textbf{not-a-knot} conditions at the two penultimate knots:
\[
s^{\left(  3\right)  }\left(  R_{1}-0\right)  =s^{\left(  3\right)  }\left(
R_{1}+0\right)  \quad\text{and }\quad s^{\left(  3\right)  }\left(
R_{N_{1}-1}-0\right)  =s^{\left(  3\right)  }\left(  R_{N_{1}-1}+0\right)  .
\]
Then the error in the interpolation satisfies
\[
\left\Vert \left(  g-s\right)  ^{\left(  r\right)  }\right\Vert _{\infty}\leq
C_{r}\left\Vert g^{\left(  4\right)  }\right\Vert _{L_{\infty}\left[
a,b\right]  }\times h^{4-r}\qquad\text{for }r=0,1,2,3,
\]
where the constant $C_{r}$ is independent of the function $g$ and the size of
the grid $h.$
\end{proposition}

For the proof we refer to \cite{deBoorNotAknot} and \cite{beatson}, inequality (2.1).

\begin{remark}
The proper value of the constant $C_{0}$ above may be estimated by the some
arguments, provided in  \cite{dahlquist}, p. 422, and Theorem $4.4.8$ therein.
Thus the value as provided in Theorem $4.4.8$ \cite{dahlquist}, seems to be a
reasonable approximation:
\[
C_{0}\leq\frac{5}{384}.
\]

\end{remark}

Now we may prove the following main theorem.

\begin{theorem}
\label{ThybridError} Let the function $f\in C^{\infty}\left(  \overline{D_{R}%
}\right)  .$ Let $N_{1}\geq5.$ Then the difference
\[
I_{\left(  N,M,K,N_{1}\right)  }^{\text{spline}}\left(  f\right)  -I_{\left(
N,M,K\right)  }^{\text{poly}}\left(  f\right)
\]
between the Hybrid Polyharmonic Cubature formula (\ref{hybridCubature}) and
the Discrete Polyharmonic Cubature formula (\ref{DiscrPolyhCubature})
satisfies the following estimate:
\begin{align*}
\left\vert I_{\left(  N,M,K,N_{1}\right)  }^{\text{spline}}\left(  f\right)
-I_{\left(  N,M,K\right)  }^{\text{poly}}\left(  f\right)  \right\vert  &
\leq\sqrt{\pi}Ch^{4}\times\left\Vert \frac{\partial^{4}f\left(  re^{i\varphi
}\right)  }{\partial r^{4}}\right\Vert _{L_{\infty}\left(  \overline{D_{R}%
}\right)  }\times\sum_{k=0}^{K}\sum_{\ell=1}^{a_{k}}\int_{0}^{R}\left\vert
w_{\left(  k,\ell\right)  }\left(  r\right)  \right\vert rdr\\
&  \leq\sqrt{\pi}Ch^{4}\times\left\Vert \frac{\partial^{4}f\left(
re^{i\varphi}\right)  }{\partial r^{4}}\right\Vert _{L_{\infty}\left(
\overline{D_{R}}\right)  }\times\left\Vert w\right\Vert .
\end{align*}
Here $C>0$ is the constant $C_{0},$ independent of the function  $f$ and the
parameter  $h,$ as provided by Proposition \ref{PBeatsonDeBoor}.
\end{theorem}%

\proof
In order to find the error bound of this approximation we recall that for
every $r$ with $0\leq r\leq R$ we have by definition (\ref{DFT}),
\[
f_{\left(  k,\ell\right)  }^{\left(  M\right)  }\left(  r\right)  =\frac{2\pi
}{M}\sum_{s=1}^{M}f\left(  re^{i\frac{2\pi}{M}s}\right)  Y_{\left(
k,\ell\right)  }\left(  \frac{2\pi}{M}s\right)  .
\]
By the linearity of the spline interpolation, for every $r$ with $0\leq r\leq
R$ we obtain
\[
SPL\left[  \left\{  f_{\left(  k,\ell\right)  }^{\left(  M\right)  }\left(
R_{m}\right)  \right\}  _{m=0}^{N_{1}}\right]  \left(  r\right)  =\frac{2\pi
}{M}\sum_{s=1}^{M}SPL\left[  \left\{  f\left(  R_{m}e^{i\frac{2\pi}{M}%
s}\right)  \right\}  _{m=0}^{N_{1}}\right]  \left(  r\right)  \times
Y_{\left(  k,\ell\right)  }\left(  \frac{2\pi}{M}s\right)  ,
\]
hence,
\begin{align}
&  SPL\left[  \left\{  f_{\left(  k,\ell\right)  }^{\left(  M\right)  }\left(
R_{m}\right)  \right\}  _{m=0}^{N_{1}}\right]  \left(  r\right)  -f_{\left(
k,\ell\right)  }^{\left(  M\right)  }\left(  r\right)  \label{SPLlinear}\\
&  =\frac{2\pi}{M}\sum_{s=1}^{M}\left\{  SPL\left[  \left\{  f\left(
R_{m}e^{i\frac{2\pi}{M}s}\right)  \right\}  _{m=0}^{N_{1}}\right]  \left(
r\right)  -f\left(  re^{i\frac{2\pi}{M}s}\right)  \right\}  \times Y_{\left(
k,\ell\right)  }\left(  \frac{2\pi}{M}s\right)  .\nonumber
\end{align}

By Proposition \ref{PBeatsonDeBoor} we obtain the inequality: \
\[
\left\vert f\left(  re^{i\frac{2\pi}{M}s}\right)  -SPL\left[  \left\{
f\left(  R_{m}e^{i\frac{2\pi}{M}s}\right)  \right\}  _{m=0}^{N_{1}}\right]
\left(  r\right)  \right\vert \leq C\left\Vert \frac{\partial^{4}f\left(
re^{i\varphi}\right)  }{\partial r^{4}}\right\Vert _{L_{\infty}\left(
\overline{D_{R}}\right)  }\times h^{4}.
\]
Hence, from (\ref{SPLlinear}) we obtain the inequality \
\begin{align*}
&  \left\vert SPL\left[  \left\{  f_{\left(  k,\ell\right)  }^{\left(
M\right)  }\left(  R_{m}\right)  \right\}  _{m=0}^{N_{1}}\right]  \left(
r\right)  -f_{\left(  k,\ell\right)  }^{\left(  M\right)  }\left(  r\right)
\right\vert \\
&  \leq\frac{2\pi}{M}C\left\Vert \frac{\partial^{4}f\left(  re^{i\varphi
}\right)  }{\partial r^{4}}\right\Vert _{L_{\infty}\left(  \overline{D_{R}%
}\right)  }\times h^{4}\sum_{s=1}^{M}\left\vert Y_{\left(  k,\ell\right)
}\left(  \frac{2\pi}{M}s\right)  \right\vert \\
&  \leq\frac{2\pi}{M}C\left\Vert \frac{\partial^{4}f\left(  re^{i\varphi
}\right)  }{\partial r^{4}}\right\Vert _{L_{\infty}\left(  \overline{D_{R}%
}\right)  }\times h^{4}\frac{M}{\sqrt{\pi}}\\
&  =2\sqrt{\pi}C\left\Vert \frac{\partial^{4}f\left(  re^{i\varphi}\right)
}{\partial r^{4}}\right\Vert _{L_{\infty}\left(  \overline{D_{R}}\right)
}\times h^{4}.
\end{align*}

Finally, we obtain the final result:
\begin{align*}
&  \left\vert I_{\left(  N,M,K,N_{1}\right)  }^{\text{spline}}\left(
f\right)  -I_{\left(  N,\infty,K\right)  }^{\text{poly}}\left(  f\right)
\right\vert \\
&  \leq\left\vert \frac{1}{2}%
{\displaystyle\sum_{k=0}^{K}}
{\displaystyle\sum_{\ell=1}^{a_{k}}}
{\displaystyle\sum_{j=1}^{N}}
\lambda_{j,\left(  k,\ell\right)  }\times t_{j,\left(  k,\ell\right)
}^{-\frac{k}{2}}\times\left\{  SPL\left[  \left\{  f_{\left(  k,\ell\right)
}^{\left(  M\right)  }\left(  R_{m}\right)  \right\}  _{m=0}^{N_{1}}\right]
\left(  \sqrt{t_{j,\left(  k,\ell\right)  }}\right)  -f_{\left(
k,\ell\right)  }^{\left(  M\right)  }\left(  \sqrt{t_{j,\left(  k,\ell\right)
}}\right)  \right\}  \right\vert \\
&  \leq\frac{1}{2}%
{\displaystyle\sum_{k=0}^{K}}
{\displaystyle\sum_{\ell=1}^{a_{k}}}
{\displaystyle\sum_{j=1}^{N}}
\left\vert \lambda_{j,\left(  k,\ell\right)  }\times t_{j,\left(
k,\ell\right)  }^{-\frac{k}{2}}\right\vert \times2\sqrt{\pi}C\left\Vert
\frac{\partial^{4}f\left(  re^{i\varphi}\right)  }{\partial r^{4}}\right\Vert
_{L_{\infty}\left(  \overline{D_{R}}\right)  }\times h^{4}.
\end{align*}
By applying inequality (\ref{GNkl}) we finish the proof.%

\endproof

\begin{remark}
From the proof above we see that there is an alternative way to understand our
scheme of obtaining the Hybrid Polyharmonic Cubature formula
(\ref{hybridCubature}): we approximate the values of the function
\[
f\left(  \sqrt{t_{j,\left(  k,\ell\right)  }}e^{i\frac{2\pi}{M}s}\right)
\]
by means of the spline values $SPL\left[  \left\{  f\left(  R_{m}%
e^{i\frac{2\pi}{M}s}\right)  \right\}  _{m=0}^{N_{1}}\right]  \left(
\sqrt{t_{j,\left(  k,\ell\right)  }}\right)  ,$ and we use this approximation
to find $f_{\left(  k,\ell\right)  }^{\left(  M\right)  }\left(
\sqrt{t_{j,\left(  k,\ell\right)  }}\right)  .$ We see that numerically, the
original scheme is more space-saving since we do not need to keep the
approximations to the set of values $\left\{  f\left(  \sqrt{t_{j,\left(
k,\ell\right)  }}e^{i\frac{2\pi}{M}s}\right)  \right\}  $ -- they are
$N\times\left(  2K-1\right)  \times M$.
\end{remark}

In our experiments below we have chosen for simplicity
\begin{equation}
R_{j}=\frac{R}{N_{1}}j\qquad\text{for }j=1,2,...,N_{1},\label{R1Rj}%
\end{equation}
and we see that
\[
\left\vert I_{\left(  N,M,K,N_{1}\right)  }^{\text{spline}}\left(  f\right)
-I_{\left(  N,M,K\right)  }^{\text{poly}}\left(  f\right)  \right\vert
\leq\sqrt{\pi}C\times\frac{R^{4}}{N_{1}^{4}}\left\Vert w\right\Vert
\times\left\Vert \frac{\partial^{4}f\left(  re^{i\varphi}\right)  }{\partial
r^{4}}\right\Vert _{L_{\infty}\left(  \overline{D_{R}}\right)  }.
\]

As proved in the following corollary, the Hybrid Polyharmonic Cubature formula
(\ref{hybridCubature}) is approximately \emph{exact} for the same subspace of
polynomials of the type $r^{2s}r^{k}Y_{\left(  k,\ell\right)  }\left(
\varphi\right)  $ which was considered in
\cite{kounchevRender2015DiscretePolyharmonicCubature}.

\begin{corollary}
\label{CexactNearly} For functions $F=r^{2s+k}Y_{k,\ell}\left(  \varphi
\right)  ,$ with $0\leq s\leq2N-1,$ $k\leq M-1-K,$ and $\ell=1,...,a_{k},$ and
for $R_{j}$ given by (\ref{R1Rj}) we have the estimate
\[
\left\vert I_{\left(  N,M,K,N_{1}\right)  }^{\text{spline}}\left(  f\right)
-I_{w}\left(  f\right)  \right\vert \leq C\times\frac{R^{2s+k}}{N_{1}^{4}%
}\left\Vert w\right\Vert ,
\]
where $C$ is the constant provided by Proposition \ref{PBeatsonDeBoor}.
\end{corollary}

\begin{remark}
We see that for the weight $w^{\left(  1\right)  }$ we have, (\ref{eqweight}%
),
\[
w^{\left(  1\right)  }=\frac{1}{r}+\cos\varphi\ =\frac{\sqrt{2\pi}}%
{r}Y_{\left(  0,1\right)  }\left(  \varphi\right)  +\sqrt{\pi}Y_{\left(
1,1\right)  }\left(  \varphi\right)
\]
hence,
\[
\left\Vert w^{\left(  1\right)  }\right\Vert =\int_{0}^{1}\frac{\sqrt{2\pi}%
}{r}rdr+\int_{0}^{1}\sqrt{\pi}rdr=\sqrt{2\pi}+\frac{1}{2}\sqrt{\pi}\approx3.4.
\]

For the weight $w^{\left(  2\right)  }$ we obtain%
\[
w_{\left(  0,1\right)  }^{\left(  1\right)  }\left(  r\right)  =\frac
{2\sqrt{2}}{\sqrt{\pi}}r\quad\text{ and }\quad w_{\left(  2k,1\right)
}^{\left(  1\right)  }\left(  r\right)  =-\frac{4}{\sqrt{\pi}}\frac{1}%
{4k^{2}-1}r\quad\text{ for }k\geq1.
\]
Hence,
\begin{align*}
\left\Vert w^{\left(  2\right)  }\right\Vert  &  =\int_{0}^{1}\frac{2\sqrt{2}%
}{\sqrt{\pi}}r^{2}dr+\sum_{k=1}^{\infty}\int_{0}^{1}\frac{4}{\sqrt{\pi}}%
\frac{1}{4k^{2}-1}r^{2}dr\\
&  =\frac{2\sqrt{2}}{3\sqrt{\pi}}+\frac{4}{3\sqrt{\pi}}\sum_{k=1}^{\infty
}\frac{1}{4k^{2}-1}\leq\frac{2\sqrt{2}}{3\sqrt{\pi}}+\frac{4}{3\sqrt{\pi}%
}\left(  \int_{2}^{\infty}\frac{1}{4x^{2}-1}dx+\frac{1}{3}\right)  \\
&  <\frac{2\sqrt{2}}{3\sqrt{\pi}}+\frac{4}{3\sqrt{\pi}}\times0.28\approx0.742
\end{align*}

\end{remark}

Corollary \ref{CexactNearly} shows that for $R=1$, $N_{1}\approx10,$ and
$C\approx\frac{5}{384}$ (by Theorem $4.4.8$ in \cite{dahlquist}),  for the
weigths $w^{\left(  1\right)  }$, $w^{\left(  2\right)  },$ the formula
$I_{\left(  N,M,K,N_{1}\right)  }^{\text{spline}}\left(  f\right)  $ is nearly
exact with precision approximately given by
\[
C\times\frac{R^{2s+k}}{N_{1}^{4}}\left\Vert w\right\Vert \leq\frac{5}%
{384}\frac{1}{10000}3.4\approx4.4271\times10^{-6}.
\]
Hence, we may consider that the coefficient of efficiency of the cubature
formula $I_{\left(  N,M,K,N_{1}\right)  }^{\text{spline}}\left(  f\right)  $
is \emph{approximately} given by
\[
E=\frac{2N\times2\left(  M-1-K\right)  }{3N_{1}\times M}.
\]
We see that if we choose $N_{1}\approx N,$ and the Fourier approximation
parameter $K$ is given, then in order to make $E\approx1$ we have to choose
\[
M\approx4\left(  K+1\right)
\]
However, unlike the references \cite{vioreanuRokhlin}, \cite{xiaoGimbutas}, we
have to take care of the error bounds which depend on the derivatives of the
functions $f$ and the asymptotic properties of the weight functions
$w_{\left(  k,\ell\right)  },$ hence we have to make even a more careful
choice of the parameters $M$ and $K.$

\section{Experimental results for the weight function $w^{\left(  1\right)
}\left(  x,y\right)  $ \label{Sexperiments}}

As in \cite{kounchevRender2015DiscretePolyharmonicCubature}, we test the
Hybrid Polyharmonic Cubature formula for integrals of the type
\[
I_{w}\left(  f\right)  =\int_{0}^{2\pi}\int_{0}^{R}f\left(  r\cos\varphi
,r\sin\varphi\right)  \cdot w^{\left(  1\right)  }\left(  r\cos\varphi
,r\sin\varphi\right)  \cdot rdrd\varphi
\]
first for the weight function
\[
w^{\left(  1\right)  }\left(  x,y\right)  =\frac{1+x}{\sqrt{x^{2}+y^{2}}%
}\text{ }=\frac{1}{r}+\cos\varphi\ =\frac{\sqrt{2\pi}}{r}Y_{\left(
0,1\right)  }\left(  \varphi\right)  +\sqrt{\pi}Y_{\left(  1,1\right)
}\left(  \varphi\right)  .
\]

We consider for our experiments three test functions:
\begin{align*}
f_{0}\left(  x,y\right)   &  =1+x^{4}+y^{3}\\
f_{1}\left(  x,y\right)   &  =1+\frac{x^{3}}{\sqrt{x^{2}+y^{2}}}+\frac{y^{7}%
}{x^{2}+y^{2}}=1+r^{2}\cos^{3}\varphi+r^{5}\sin^{7}\varphi\\
f_{2}\left(  x,y\right)   &  =\cos\left(  ax+by\right)  \text{ for }a=10\text{
and }b=20.
\end{align*}
The first test function $f_{0}$ is a polynomial of degree $4$, the second
$f_{1}$ is not smooth at $0$ and the last one $f_{2}$ is of oscillatory type. 

\subsection{Results for the Hybrid Polyharmonic Cubature Formula}

As in \cite{kounchevRender2015DiscretePolyharmonicCubature}, the number $M$ of
points on the circles is choosen to be equal to $9,25,63,83$ and the number
$N$ of concentric circles is chosen to be equal to $10,15,25,35,50.$ The
Gauss-Jacobi quadrature rules are with $N_{1}=N$ points. We have the exact
values given by
\begin{align*}
I_{1}\left(  f_{0}w^{\left(  1\right)  }\right)    & =\frac{43}{20}\pi
\approx6.754424205218060\\
I_{1}\left(  f_{1}w^{\left(  1\right)  }\right)    & =\frac{35}{16}\pi
\approx6.872\,233\,929727\,67\\
I_{1}\left(  f_{2}w^{\left(  1\right)  }\right)    & \approx
0.301\,310\,995335\,215\,
\end{align*}

Below is the table with the values of $I_{\left(  N,M,K,N_{1}\right)
}^{\text{spline}}\left(  f_{0}\right)  ,$ and the error.  We see that the
error is bigger compared with the similar table in
\cite{kounchevRender2015DiscretePolyharmonicCubature}.
\begin{gather*}%
\begin{tabular}
[c]{|lllll|}\hline
N/M & \multicolumn{1}{|l}{\textbf{9}} & \multicolumn{1}{|l}{\textbf{25}} &
\multicolumn{1}{|l}{\textbf{63}} & \multicolumn{1}{|l|}{\textbf{83}}\\\hline
\textbf{10} & \multicolumn{1}{|r|}{6,754363639426710} &
\multicolumn{1}{c|}{6,754363639426710} &
\multicolumn{1}{c|}{6,754363639426710} &
\multicolumn{1}{|l|}{6,754363639426710}\\\hline
\textbf{15} & \multicolumn{1}{c|}{6,754415757033810} &
\multicolumn{1}{c|}{6,754415757033810} &
\multicolumn{1}{c|}{6,754415757033810} &
\multicolumn{1}{|l|}{6,754415757033800}\\\hline
\textbf{25} & \multicolumn{1}{|c|}{6,754423468244570} &
\multicolumn{1}{c|}{6,754423468244570} &
\multicolumn{1}{c|}{6,754423468244570} &
\multicolumn{1}{r|}{6,754423468244570}\\\hline
\textbf{35} & \multicolumn{1}{c|}{6,75424054924440} &
\multicolumn{1}{c|}{6,754424054924440} &
\multicolumn{1}{c|}{6,754424054924430} &
\multicolumn{1}{|l|}{6,754424054924440}\\\hline
\textbf{50} & \multicolumn{1}{|c|}{6,754424177151970} &
\multicolumn{1}{c|}{6,754424177151970} &
\multicolumn{1}{c|}{6,754424177151970} &
\multicolumn{1}{r|}{6,754424177151970}\\\hline
Error & \multicolumn{1}{|l}{} & \multicolumn{1}{|l}{} & \multicolumn{1}{|l}{}
& \multicolumn{1}{|l|}{}\\\hline
N/M & \multicolumn{1}{|l}{\textbf{9}} & \multicolumn{1}{|l}{\textbf{25}} &
\multicolumn{1}{|l}{\textbf{63}} & \multicolumn{1}{|l|}{\textbf{83}}\\\hline
\textbf{10} & \multicolumn{1}{|c|}{-0,000060565791338} &
\multicolumn{1}{c|}{-0,000060565791338} &
\multicolumn{1}{c|}{-0,000060565791338} &
\multicolumn{1}{r|}{-0,000060565791338}\\\hline
\textbf{15} & \multicolumn{1}{|c|}{-0,000008448184238} &
\multicolumn{1}{c|}{-0,000008448184238} &
\multicolumn{1}{c|}{-0,000008448184238} &
\multicolumn{1}{r|}{-0,000008448184248}\\\hline
\textbf{25} & \multicolumn{1}{|c|}{-0,000000736973478} &
\multicolumn{1}{c|}{-0,000000736973478} &
\multicolumn{1}{c|}{-0,000000736973478} &
\multicolumn{1}{r|}{-0,000000736973478}\\\hline
\textbf{35} & \multicolumn{1}{|c|}{-0,000000150293608} &
\multicolumn{1}{c|}{-0,000000150293608} &
\multicolumn{1}{c|}{-0,000000150293618} &
\multicolumn{1}{r|}{-0,000000150293608}\\\hline
\textbf{50} & \multicolumn{1}{|c|}{-0,000000028066078} &
\multicolumn{1}{c|}{-0,000000028066078} &
\multicolumn{1}{c|}{-0,000000028066078} &
\multicolumn{1}{r|}{-0,000000028066078}\\\hline
\end{tabular}
\\
\text{\textbf{Hybrid Polyharmonic cubature} for }f_{0}w^{\left(  1\right)
}=\left(  1+x^{4}+y^{3}\right)  \left(  \frac{1}{r}+\cos\varphi\right)  \text{
and the error.}\\
\text{\textbf{True value} is }\frac{43}{20}\pi\approx6.754424205218060
\end{gather*}

Below is the table with the values of $I_{\left(  N,M,K,N_{1}\right)
}^{\text{spline}}\left(  f_{1}\right)  $ and the error. Again,  the error is
bigger compared with the similar table in
\cite{kounchevRender2015DiscretePolyharmonicCubature}:%

\begin{gather*}%
\begin{tabular}
[c]{|l|l|l|l|l|}\hline
N/M & \textbf{9} & \textbf{25} & \textbf{63} & \textbf{83}\\\hline
\textbf{10} & 6,87224296287783 & 6,87224296287783 & 6,87224296287783 &
6,87224296287783\\\hline
\textbf{15} & 6,87223588060173 & 6,87223588060173 & 6,87223588060173 &
6,87223588060173\\\hline
\textbf{25} & 6,87223420205342 & 6,87223420205342 & 6,87223420205342 &
6,87223420205342\\\hline
\textbf{35} & 6,87223400297000 & 6,87223400297000 & 6,87223400297000 &
6,87223400297000\\\hline
\textbf{50} & 6,87223394775545 & 6,87223394775545 & 6,87223394775545 &
6,87223394775545\\\hline
Error & \multicolumn{1}{|l|}{} & \multicolumn{1}{|l|}{} &
\multicolumn{1}{|l|}{} & \multicolumn{1}{|l|}{}\\\hline
N/M & \multicolumn{1}{|l|}{\textbf{9}} & \multicolumn{1}{|l|}{\textbf{25}} &
\multicolumn{1}{|l|}{\textbf{63}} & \multicolumn{1}{|l|}{\textbf{83}}\\\hline
\textbf{10} & \multicolumn{1}{|c|}{0,00000903315016} &
\multicolumn{1}{c|}{0,00000903315016} & \multicolumn{1}{c|}{0,00000903315016}
& \multicolumn{1}{r|}{0,00000903315016}\\\hline
\textbf{15} & \multicolumn{1}{|c|}{0,00000195087406} &
\multicolumn{1}{c|}{0,00000195087406} & \multicolumn{1}{c|}{0,00000195087406}
& \multicolumn{1}{r|}{0,00000195087406}\\\hline
\textbf{25} & \multicolumn{1}{|c|}{0,00000027232575} &
\multicolumn{1}{c|}{0,00000027232575} & \multicolumn{1}{c|}{0,00000027232575}
& \multicolumn{1}{r|}{0,00000027232575}\\\hline
\textbf{35} & \multicolumn{1}{|c|}{0,00000007324233} &
\multicolumn{1}{c|}{0,00000007324233} & \multicolumn{1}{c|}{0,00000007324233}
& \multicolumn{1}{r|}{0,00000007324233}\\\hline
\textbf{50} & \multicolumn{1}{|c|}{0,00000001802778} &
\multicolumn{1}{c|}{,00000001802778} & \multicolumn{1}{c|}{0,00000001802778} &
\multicolumn{1}{r|}{0,00000001802778}\\\hline
\end{tabular}
\\
\text{\textbf{Hybrid Polyharmonic cubature} for }f_{1}w^{\left(  1\right)
}=\left(  1+r^{2}\cos^{3}\varphi+r^{5}\sin^{7}\varphi\right)  \left(  \frac
{1}{r}+\cos\varphi\right)  \\
\text{\textbf{True value} is }\frac{35}{16}\pi\approx
6.\,872\,233\,929\,727\,67
\end{gather*}

There is no mistake in the table above since results do not change with $M$!
This is due to the fact that%
\[
\sin^{7}\varphi=\frac{7}{64}\sin5\varphi-\frac{21}{64}\sin3\varphi-\frac
{1}{64}\sin7\varphi+\frac{35}{64}\sin\varphi
\]
which implies that the function $g=$ $f_{1}w^{\left(  1\right)  }$ has no
terms $g_{\left(  k,\ell\right)  }\left(  r\right)  $ with $k\geq4;$ hence,
the interpolation cubic splines of the polynomials $g_{\left(  k,\ell\right)
}\left(  r\right)  r$ coincide with them.

For the oscillatory test function $f_{2}$ we obtain%

\begin{gather*}%
\begin{tabular}
[c]{ccccc}\hline
N%
$\backslash$%
M & \textbf{9} & \textbf{25} & \textbf{63} & \multicolumn{1}{|l}{\textbf{83}%
}\\\hline
\multicolumn{1}{|l}{\textbf{10}} & \multicolumn{1}{|c|}{0,19210087475239} &
\multicolumn{1}{c|}{0,58134400368821} & \multicolumn{1}{c|}{0,56846433865624}
& 0,56846433865624\\\hline
\multicolumn{1}{|l}{\textbf{15}} & \multicolumn{1}{c|}{-0,00842490123728} &
\multicolumn{1}{c|}{0,38518390970894} & \multicolumn{1}{c|}{0,37239275649383}
& 0,37239275649383\\\hline
\multicolumn{1}{|l}{\textbf{25}} & \multicolumn{1}{|c|}{-0,08069670830424} &
\multicolumn{1}{c|}{0,31431334300881} & \multicolumn{1}{c|}{0,30152604401835}
& \multicolumn{1}{c|}{0,30152604401835}\\\hline
\multicolumn{1}{|l}{\textbf{35}} & \multicolumn{1}{c|}{-0,08150067781664} &
\multicolumn{1}{c|}{0,31360790761542} & \multicolumn{1}{c|}{0,30081999553130}
& 0,30081999553130\\\hline
\multicolumn{1}{|l}{\textbf{50}} & \multicolumn{1}{|c|}{-0,08116599113863} &
\multicolumn{1}{c|}{0,31395572503057} & \multicolumn{1}{c|}{0,30116759220177}
& \multicolumn{1}{c|}{0,30116759220177}\\\hline
Error &  &  &  & \\\hline
\multicolumn{1}{|l}{N%
$\backslash$%
M} & \textbf{9} & \textbf{25} & \textbf{63} & \multicolumn{1}{|l}{\textbf{83}%
}\\\hline
\multicolumn{1}{|l}{\textbf{10}} & \multicolumn{1}{|c|}{-0,10921012058282} &
\multicolumn{1}{c|}{0,28003300835300} & \multicolumn{1}{c|}{0,26715334332102}
& \multicolumn{1}{c|}{0,26715334332102}\\\hline
\multicolumn{1}{|l}{\textbf{15}} & \multicolumn{1}{|c|}{-0,30973589657249} &
\multicolumn{1}{c|}{0,08387291437372} & \multicolumn{1}{c|}{0,07108176115862}
& \multicolumn{1}{c|}{0,07108176115862}\\\hline
\multicolumn{1}{|l}{\textbf{25}} & \multicolumn{1}{|c|}{-0,38200770363946} &
\multicolumn{1}{c|}{0,01300234767360} & \multicolumn{1}{c|}{0,00021504868313}
& \multicolumn{1}{c|}{0,00021504868313}\\\hline
\multicolumn{1}{|l}{\textbf{35}} & \multicolumn{1}{|c|}{-0,38281167315186} &
\multicolumn{1}{c|}{0,01229691228021} & \multicolumn{1}{c|}{-0,00049099980392}
& \multicolumn{1}{c|}{-0,00049099980392}\\\hline
\multicolumn{1}{|l}{\textbf{50}} & \multicolumn{1}{|c|}{-0,38247698647384} &
\multicolumn{1}{c|}{0,01264472969535} & \multicolumn{1}{c|}{-0,00014340313345}
& \multicolumn{1}{c|}{-0,00014340313345}\\\hline
\end{tabular}
\\
\text{\textbf{Hybrid Polyharmonic cubature} for }f_{2}w^{\left(  1\right)
}=\left(  \cos\left(  10x+20y\right)  \right)  \left(  \frac{1}{r}+\cos
\varphi\right) \\
\text{\textbf{True value} is }I_{1}\left(  f_{2}w^{\left(  1\right)  }\right)
\approx0.301\,310\,995\,335\,215
\end{gather*}

\section{Experimental results for the weight function $w^{\left(  2\right)
}\left(  x,y\right)  $ \label{Sexperiments2}}

As in \cite{kounchevRender2015DiscretePolyharmonicCubature}, we consider the
second weight function:
\begin{equation}
w^{\left(  2\right)  }\left(  re^{i\varphi}\right)  :=\left\vert y\right\vert
=\left\vert r\sin\varphi\right\vert \label{w2}%
\end{equation}
We have its orthonormalized Fourier coefficients:
\[
w_{\left(  0,1\right)  }\left(  r\right)  =\frac{2\sqrt{2}}{\sqrt{\pi}}r\text{
and }w_{\left(  2k,1\right)  }\left(  r\right)  =-\frac{4}{\sqrt{\pi}}\frac
{1}{4k^{2}-1}r\text{ for }k\geq1.
\]

\section{Results for the Hybrid Discrete Polyharmonic Cubature}

As in \cite{kounchevRender2015DiscretePolyharmonicCubature}, we consider the
test functions
\[
f_{3}\left(  x,y\right)  =30x^{12},\qquad f_{4}\left(  x,y\right)  =\left\vert
y\right\vert ,
\]
for which
\begin{align*}
I_{w^{\left(  2\right)  }}\left(  f\right)   &  =I_{1}\left(  30x^{12}%
\left\vert y\right\vert \right)  =\frac{8}{13}\approx
0.615\,384\,615\,384\,616\,\\
I_{1}\left(  \left\vert y\right\vert \left\vert y\right\vert \right)   &
=\int_{0}^{1}\int_{0}^{2\pi}\left\vert r\sin\varphi\right\vert ^{2}d\varphi
rdr=\frac{1}{4}\pi\approx0.785\,398\,163\,397\,448\,
\end{align*}

The table contains the results for the \textbf{Hybrid} \textbf{Discrete
Polyharmonic Cubature }for the test function $f_{3}\left(  x,y\right)  $
\textbf{ }(\ref{hybridCubature}), where $K=22$\textbf{:\ }%
\begin{gather*}%
\begin{tabular}
[c]{ccccc}\hline
N%
$\backslash$%
M & \multicolumn{1}{|c|}{\textbf{9}} & \multicolumn{1}{c|}{\textbf{25}} &
\multicolumn{1}{c|}{\textbf{63}} & \textbf{83}\\\hline
\multicolumn{1}{|c}{\textbf{10}} & \multicolumn{1}{c|}{0,565617343585166} &
\multicolumn{1}{c|}{0,620572422003199} &
\multicolumn{1}{c|}{0,620572422003199} & 0,620572422003199\\\hline
\multicolumn{1}{|c}{\textbf{15}} & \multicolumn{1}{|c|}{0,561709413462587} &
\multicolumn{1}{c|}{0,616243839415133} &
\multicolumn{1}{c|}{0,616243839415133} &
\multicolumn{1}{c|}{0,616243839415132}\\\hline
\multicolumn{1}{|c}{\textbf{25}} & \multicolumn{1}{c|}{0,561006480042405} &
\multicolumn{1}{c|}{0,615463257008362} &
\multicolumn{1}{c|}{0,615463257008361} & 0,615463257008362\\\hline
\multicolumn{1}{|c}{\textbf{35}} & \multicolumn{1}{|c|}{0,560949764331337} &
\multicolumn{1}{c|}{0,615400379071044} &
\multicolumn{1}{c|}{0,615400379071044} &
\multicolumn{1}{c|}{0,615400379071044}\\\hline
\multicolumn{1}{|c}{\textbf{50}} & \multicolumn{1}{|c|}{0,560937835197964} &
\multicolumn{1}{c|}{0,615387283068315} &
\multicolumn{1}{c|}{0,615387283068315} &
\multicolumn{1}{c|}{0,615387283068316}\\\hline
Error &  &  &  & \\\hline
N%
$\backslash$%
M & \multicolumn{1}{|c|}{\textbf{9}} & \multicolumn{1}{c|}{\textbf{25}} &
\multicolumn{1}{c|}{\textbf{63}} & \textbf{83}\\\hline
\multicolumn{1}{|c}{\textbf{10}} & \multicolumn{1}{|c|}{-0,0497672717994500} &
\multicolumn{1}{c|}{0,0051878066185830} &
\multicolumn{1}{c|}{0,0051878066185830} &
\multicolumn{1}{c|}{0,0051878066185830}\\\hline
\multicolumn{1}{|c}{\textbf{15}} & \multicolumn{1}{|c|}{-0,0536752019220290} &
\multicolumn{1}{c|}{0,0008592240305171} &
\multicolumn{1}{c|}{0,0008592240305171} &
\multicolumn{1}{c|}{0,0008592240305161}\\\hline
\multicolumn{1}{|c}{\textbf{25}} & \multicolumn{1}{|c|}{-0,0543781353422109} &
\multicolumn{1}{c|}{0,0000786416237460} &
\multicolumn{1}{c|}{0,0000786416237450} &
\multicolumn{1}{c|}{0,0000786416237460}\\\hline
\multicolumn{1}{|c}{\textbf{35}} & \multicolumn{1}{|c|}{-0,0544348510532789} &
\multicolumn{1}{c|}{0,0000157636864280} &
\multicolumn{1}{c|}{0,0000157636864280} &
\multicolumn{1}{c|}{0,0000157636864280}\\\hline
\multicolumn{1}{|c}{\textbf{50}} & \multicolumn{1}{|c|}{-0,0544467801866519} &
\multicolumn{1}{c|}{0,0000026676836991} &
\multicolumn{1}{c|}{0,0000026676836991} &
\multicolumn{1}{c|}{0,0000026676837001}\\\hline
\end{tabular}
\\
\text{\textbf{Hybrid Discrete Polyharmonic cubature} for }30x^{12}\left\vert
y\right\vert \\
\text{\textbf{True value} is }\approx0.615\,384\,615\,\allowbreak384\,616
\end{gather*}

For the test function $f_{4}\left(  x,y\right)  =\left\vert y\right\vert ,$
with $K=22$ in the formula\textbf{ }(\ref{hybridCubature}), we obtain the table

\begin{gather*}%
\begin{tabular}
[c]{rrrrc}\hline
N%
$\backslash$%
M & \multicolumn{1}{|r|}{\textbf{9}} & \multicolumn{1}{r|}{\textbf{25}} &
\multicolumn{1}{r|}{\textbf{63}} & \textbf{83}\\\hline
\multicolumn{1}{|r}{\textbf{10}} & \multicolumn{1}{r|}{0,785206660} &
\multicolumn{1}{c}{0,785352337} & \multicolumn{1}{c}{0,785367124} &
0,785369362\\\hline
\multicolumn{1}{|r}{\textbf{15}} & \multicolumn{1}{|r|}{0,785208297} &
\multicolumn{1}{c}{0,78535897} & \multicolumn{1}{c}{0,785373081} &
\multicolumn{1}{c|}{0,785375274}\\\hline
\multicolumn{1}{|r}{\textbf{25}} & \multicolumn{1}{r|}{0,785208235} &
\multicolumn{1}{c}{0,785361119} & \multicolumn{1}{c}{0,785374994} &
0,785377171\\\hline
\multicolumn{1}{|r}{\textbf{35}} & \multicolumn{1}{|r|}{0,785208149} &
\multicolumn{1}{c}{0,78536144} & \multicolumn{1}{c}{0,785375276} &
\multicolumn{1}{c|}{0,785377452}\\\hline
\multicolumn{1}{|r}{\textbf{50}} & \multicolumn{1}{|r|}{0,785208109} &
\multicolumn{1}{c}{0,785361541} & \multicolumn{1}{c}{0,785375364} &
\multicolumn{1}{c|}{0,785377539}\\\hline
Error &  &  &  & \\\hline
N%
$\backslash$%
M & \multicolumn{1}{|c|}{\textbf{9}} & \multicolumn{1}{c|}{\textbf{25}} &
\multicolumn{1}{c|}{\textbf{63}} & \textbf{83}\\\hline
\multicolumn{1}{|c}{\textbf{10}} & \multicolumn{1}{|c|}{-0,0001915037} &
\multicolumn{1}{c|}{-0,0000458267} & \multicolumn{1}{c|}{-0,0000310393} &
\multicolumn{1}{c|}{-0,0000288009}\\\hline
\multicolumn{1}{|c}{\textbf{15}} & \multicolumn{1}{|c|}{-0,0001898666} &
\multicolumn{1}{c|}{-0,0000391939} & \multicolumn{1}{c|}{-0,0000250824} &
\multicolumn{1}{c|}{-0,0000228890}\\\hline
\multicolumn{1}{|c}{\textbf{25}} & \multicolumn{1}{|c|}{-0,0001899281} &
\multicolumn{1}{c|}{-0,0000370443} & \multicolumn{1}{c|}{-0,0000231697} &
\multicolumn{1}{c|}{-0,0000209919}\\\hline
\multicolumn{1}{|c}{\textbf{35}} & \multicolumn{1}{|c|}{-0,0001900142} &
\multicolumn{1}{c|}{-0,0000367231} & \multicolumn{1}{c|}{-0,0000228870} &
\multicolumn{1}{c|}{-0,0000207116}\\\hline
\multicolumn{1}{|c}{\textbf{50}} & \multicolumn{1}{|c|}{-0,0001900544} &
\multicolumn{1}{c|}{-0,0000366227} & \multicolumn{1}{c|}{-0,0000227993} &
\multicolumn{1}{c|}{-0,0000206248}\\\hline
\end{tabular}
\\
\text{\textbf{Hybrid Polyharmonic cubature} for }=\left\vert y\right\vert
^{2}\\
\text{\textbf{True value} is }\approx\allowbreak0.785\,398\,163\,\allowbreak
397\,448
\end{gather*}

\section{Conclusions}

The comparison with the results in
\cite{kounchevRender2015DiscretePolyharmonicCubature} shows that the Hybrid
Polyharmonic Cubature Formula provides a satisfactory tradeoff between
precision and number of evaluation points.

\begin{acknowledgement}
Both authors acknowledge the partial support by the Bulgarian NSF
Grant\textbf{ I02/19, 2015}. The first named author acknowledges partial
support by the Humboldt Foundation. 
\end{acknowledgement}

\end{document}